\documentclass[11pt]{article}

\usepackage[english]{babel}
\usepackage{amsmath,amssymb,amsfonts}
\newtheorem{theorem}{Theorem}[section]
\newtheorem{lemma}[theorem]{Lemma}
\newtheorem{proposition}[theorem]{Proposition}
\newtheorem{definition}{Definition}[section]
\newtheorem{example}[theorem]{Example}

\newtheorem{remark}[theorem]{Remark}
\newtheorem{corollary}[theorem]{Corollary}

\def\Dim{\noindent\hbox{{\bf Proof.}$\;\; $}}          
\def\finedim{{\hfill\hbox{\enspace${ \square}$}} \smallskip}    
\def\sqr#1#2{{\vcenter{\vbox{\hrule height .#2pt
     \hbox{\vrule width .#2pt height#1pt \kern#1pt \vrule
     width .#2pt} \hrule height .#2pt}}}}
\def\square{\mathchoice\sqr54\sqr54\sqr{4.1}3\sqr{3.5}3}

\makeatletter
\@addtoreset{equation}{section}

\makeatother

\def\sqr#1#2{{\vcenter{\vbox{\hrule height .#2pt \hbox{\vrule
 width .#2pt height#1pt \kern#1pt \vrule
width .#2pt} \hrule height .#2pt}}}}
\def\square{\mathchoice\sqr54\sqr54\sqr{4.1}3\sqr{3.5}3}

\def\ds{\begin{displaystyle}}
\def\eds{\end{displaystyle}}

\def\<{\langle }
\def\>{\rangle }

\def\E{\mathbb E}

\def\cald{{\cal D}}

\def\calf{{\cal F}}

\def\call{{\cal L}}
\def\cals{{\cal S}}

\title{Stochastic LQ and Associated Riccati equation of  PDEs Driven by State- and Control-Dependent White Noise}

\author{Ying Hu
	\thanks{
		IRMAR, Universit{\'e} Rennes 1, Campus de Beaulieu, 35042 Rennes
		Cedex, France, and School of Mathematical Sciences, Fudan
		University, Shanghai 200433, China. Partially supported by Lebesgue
		Center of Mathematics ``Investissements d'avenir"
		program-ANR-11-LABX-0020-01, by ANR CAESARS (Grant No. 15-CE05-0024)
		and by ANR MFG (Grant No. 16-CE40-0015-01). email:
		\texttt{ying.hu@univ-rennes1.fr} } \and
	Shanjian Tang
	\thanks{
		Department of Finance and Control Sciences, School of Mathematical
		Sciences, Fudan University, Shanghai 200433, China. Partially
		supported by National Science Foundation of China (Grant No.
		11631004) and Science and Technology Commission of Shanghai
		Municipality (Grant No. 14XD1400400). email:
		\texttt{sjtang@fudan.edu.cn}}}

\begin{document}
\maketitle

\begin{abstract}
	The optimal stochastic control problem with a quadratic cost functional for  linear partial differential equations (PDEs) driven by a state- and control-dependent white noise is formulated and studied.  Both finite- and infinite-time horizons are considered. The multiplicative white noise dynamics of the system give rise  to a new phenomenon of  singularity to the associated Riccati equation and even to the Lyapunov equation. Well-posedness of both Riccati equation and Lyapunov equation are obtained  for the first time. The linear feedback coefficient of the optimal  control turns out to be singular and expressed  in terms of the solution of the associated Riccati equation. The null controllability is shown to be equivalent to the existence of the solution to Riccati equation with the singular terminal value.  Finally, the controlled Anderson model is addressed as an illustrating example.  
	
	
	
\noindent\textit{Short title}: LQ control of white noise-driven PDEs. 
	
\end{abstract}

\section{Introduction}

In this paper, we consider the following stochastic evolutionary equation driven by both state- and control-dependent white noise:
\begin{equation*}
dX_t=(AX_t+B_tu_t)\,dt+\sum_{j=1}^\infty (C_j(t)X_t+D_j(t)u_t)\,d \beta_t^j,\quad X_0=x\in H
\end{equation*}
where 
$A$ is the infinitesimal generator of a strongly continuous
semigroup $e^{tA}$ of linear operators,   $B, C_j, $ and $D_j$ are some bounded operators, and  $W$ is a cylindrical
Wiener process in a Hilbert space $H$, with $\{\beta^j(t):=\langle W(t), e_j\rangle, j=1,2,\ldots\}$ being independent Brownian motions for an orthonormal basis $\{e_j, j=1,2,\ldots\}$ of $H$. 
The cost functional is
\begin{equation*}
J(x,u)=\E\int_0^T[\langle Q_tX_t,X_t\rangle+\langle R_tu_t,u_t\rangle]\,dt+\E[\langle GX_T,X_T\rangle],
\end{equation*}
where $Q, G$, and $R$ are some bounded operators. The optimal control problem is to find a $U$-valued  adapted square-integrable process  ${\overline u}$ in a feedback form (via the associated Riccati equation) such that $J(x, {\overline u})$ is the minimal value of the cost functional $J(x, \cdot)$. More precise formulation will be given in the next section. 

The general  theory of linear quadratic optimal  control  (the so-called LQ theory) of Kalman~\cite{Kalman} paved one mile stone in the deterministic optimal control theory. The general stochastic extension in a Euclidean space was given by Wonham~\cite{Wonham1968} for the deterministic coefficients,  and  was further developed by Bismut~\cite{Bismut1976} for the random coefficients. Subsequently, it was further studied by Peng~\cite{Peng1992} and Tang~\cite{Tang2003},  and its theory is now rather complete.  

 Ichikawa~\cite{Ichikawa1976,Ichikawa1979} considered  the infinite-dimensional  extension of Kalman's LQ theory under the following setting:  H is an infinite-dimensional Hilbert space and C is a bounded linear operator. Da Prato and Ichikawa~\cite{DaIchikawa1985} studied the infinite-dimensional LQ  problem for the case of  $D=0$,  self-adjoint $A$, and unbounded coefficient $B$. The infinite dimensional case with  stochastic coefficients driven by the so-called colored noise (where $C$ is a Hilbert-Schmidt operator) is referred to Guatteri and Tessitore~\cite{GuatteriTessitore}. To our best knowledge, all the above-mentioned papers are restricted within the case when the linear SPDEs are driven by the so-called colored noise, which excludes the celebrated Anderson model. In this paper, we address the infinite dimensional stochastic LQ problem driven by an infinite number of Brownian motions (the so-called space-time white noise).  
 

The introduction of the space-time white noise leads to the difficulty  that  the infinite sum $\sum_{i=1}^\infty C_i^*(s)P_s C_i(s)$ appears in both  associated Lyapunov equation \eqref{lyapunov} and Riccati equation~\eqref{Riccati 2},  and thus challenges the solvability of both equations. To overcome this difficulty for Lyapunov equation \eqref{lyapunov}, we introduce the representation via the solution of forward SPDE to establish an estimate of the sum, and for more details,  see our Proposition~\ref{Prop3.3} and its proof. It is conventional to study the Riccati equation via the quasi-linearization method. While in our context of the space-time white noise,  the coefficients of these quasi-linearized equations become singular in the sense that these coefficients explode at both ends (time $0$ and time $T$). Some fine estimates are applied to deduce the monotonicity and convergence of solutions of quasi-linearized equations. For more details, see our Theorem~\ref{Th4.4} and its proof.   Finally, due to the space-time white noise in our context, the conventional Yosida's approximation could not be applied to get the energy equality, and to attack the new difficulty,  a new truncation is carefully constructed to deduce the energy equality and thus the feedback law of the optimal control. For more details, see our Theorem~\ref{EnergyEquality} and its proof. 

We note that  Anderson model has been widely studied in the litterature, and for more details, see Carmona and Molchanov~\cite{CarmonaMalchanov1994},  Conus, Joseph, and  Khoshnevisan~\cite{ConusJosephKhoshnevisan2013}, and the references therein. We also emphasize that our results succeed at  inclusion of the controlled Anderson SPDE. See Section 8. 

The paper is organized as follows.  In Section 2, we give the precise formulation of our quadratic optimal stochastic control problem for  linear partial differential equations driven by a white noise. In Section 3, we study well-posedness of Lyapunov equations. In Section 4, we study the associated Riccati equation. In Section 5, we characterize the optimal control as a feedback form via the solution of Riccati equation. In Section 6, we address the infinite-horizon LQ control problem for the case of time-invariant coefficients. We show that when the system is stabilizable, the associated algebraic Riccati equation has a unique solution, and is again used to synthesize the optimal control into a feedback form.  In Section 7, the null controllability is proved to be equivalent to the existence of solution of Riccati equation with the singular terminal condition. Finally in Section 8, we give  examples for the controlled Anderson model. 

\section{Formulation of the linear quadratic optimal control}

Let $H,U$ be two  separable Hilbert spaces. By $\cals (H)$,  we denote the space of all self-adjoint and  bounded linear operators on H，and by $\cals^+(H)$ we denote the set of all non-negative operators in  $\cals (H)$. Moreover, if $I\subset \mathbb{R}^+$  is an interval (bounded or unbounded), we denote by $C_s(I;\cals(H))$ (resp. $C_s(I;\cals^+(H))$)  the set of all maps $f: I\to \cals(H)$  (resp. $f: I\to \cals^+(H)$) such that $f(\cdot)$  is strongly continuous in $H$.

Consider the following stochastic evolutionary equation driven by both state- and control-dependent white noise:
\begin{equation}\label{control system}
dX_t=(AX_t+B_tu_t)\,dt+\sum_{j=1}^\infty (C_j(t)X_t+D_j(t)u_t)\,d \beta_t^j,\quad X_0=x\in H, 
\end{equation}
which has the following mild form:
\begin{equation}\label{eq:state mild}
X_t=e^{At}x+\int_0^t e^{A(t-s)}B_su_s\,ds+\int_0^t\sum_{j=1}^\infty e^{A(t-s)}  (C_j(s)X_s+D_j(s) u_s)\,d \beta_s^j.
\end{equation}
Here, 
$A$ is the infinitesimal generator of a strongly continuous
semigroup $e^{tA}$ of linear operators,  $B\in L^\infty(0,T; \call(U,H)), C_j\in  L^\infty(0,T;\call(H)), D_j\in L^\infty(0,T; \call  (U,H))$
with the standard assumption that for  some $\alpha \in (0,\frac{1}{ 2})$ and $c>0$, 
\begin{equation}\label{AC0}
\sum_{j=1}^\infty |e^{At}C_j(s)x|^2_H\le ct^{-2\alpha}|x|^2_H, \quad t>0. 
\end{equation}
$W$ is a cylindrical
Wiener process in $H$,  $\{\beta^j(t):=\langle W(t), e_j\rangle, j=1,2,\ldots\}$ are independent Brownian motions, with $\{e_j, j=1,2,\ldots\}$ being an orthonormal basis of $H$. 
The cost functional is
\begin{equation}\label{eq:cost 1}
J(x,u)=\E\int_0^T[\langle Q_tX_t,X_t\rangle+\langle R_tu_t,u_t\rangle]\,dt+\E[\langle GX_T,X_T\rangle], \quad u\in L^2_{\calf}(0,T;U)
\end{equation}
where $Q\in L^\infty(0,T; \cals^+(H)), G\in \cals^+(H)$, and $R\in  L^\infty(0,T; \cals^+(U))$ is strictly positive in the following sense: there is a positive number $\delta$ such that $R\ge \delta I_U$. Throughout the paper, we  assume that  for any $v\in U$, there is a constant $c>0$ such that
\begin{equation}\label{control}
\sum_{j=1}^\infty|D_jv|_H^2\le c|v|^2_U. 
\end{equation}

The optimal control problem is to minimize $J(x,\cdot)$ among all the  controls in $L^2_{\calf}(0,T;U)$. 

\begin{remark} Condition~\eqref{control} means that $D=(D_1, D_2, \ldots)$ is a Hilbert-Schmidt operator. It is still open how to replace this condition with a condition like~\eqref{AC0}.
\end{remark}

\begin{lemma}\label{forward}  For $u\in \call^2_{\calf}(0,T;U)$, the system~\eqref{control system} has a unique mild solution $X$ in the space $C_\calf([0,T]; L^2(\Omega, H))$ such that for some $C>0$,
		\begin{eqnarray}
	\sup_{0\le t\le T}\mathbb{E}[|X_t|^2]&\le& C \left(|x|^2+\mathbb{E}\int_0^T|u_s|^2ds \right).
	\end{eqnarray}	
	\end{lemma}

\Dim The existence and uniqueness of the mild solution can be found in \cite{DaZa}. We now derive the desired estimate for the solution.  From~\eqref{eq:state mild}, \eqref{AC0}, and \eqref{control}, we have
\begin{eqnarray*}
	\mathbb{E}[|X_t|^2]&\le& C\left(|x|^2+\mathbb{E}\int_0^t|u_s|^2ds\right)\nonumber \\
	 &&+ C \mathbb{E}\int_0^t\left(\sum_{j=1}^\infty|e^{A(t-s)}C_j(s)X_s|^2+\sum_{j=1}^\infty|e^{A(t-s)}D_j(s)u_s|^2\right)\, ds\\
	 &\le& C \left(|x|^2+\mathbb{E}\int_0^t|u_s|^2ds\right)+\nonumber  C \mathbb{E}\int_0^t\sum_{j=1}^\infty (t-s)^{-2\alpha}|X_s|^2_H\, ds.
	 \end{eqnarray*}
Using  an extended Gronwall's inequality (see, e.g. \cite{Henry}), we have the desired estimate. 
\finedim

\section{Lyapunov equation: existence and uniqueness of solutions}

We first give results on Lyapunov equation, which will be needed in the study of Riccati equation. 

\subsection{Forward SDE}

Let $A_0, {\widehat C}_j\in L^\infty(0,T; \call(H))$ with $j=1,2,\ldots$. 

Assume that for some number $c>0$, 
\begin{equation}\label{A_0}
|A_0(s)|_{\call{(H)}}\le c (T-s)^{-\alpha}, \quad s\in  [0,T), 
\end{equation}
and 
\begin{equation}\label{AC}
\sum_{j=1}^\infty |e^{At}{\widehat C}_j(s)x|^2_H\le c \left(t^{-2\alpha}+(T-s)^{-2\alpha}\right)|x|^2_H, \quad (t,s)\in (0 ,\infty)\times [0,T).
\end{equation}

\begin{remark} Assumptions~\eqref{A_0} and \eqref{AC} are introduced to study the quasi-linearized sequence of Lyapunov equations for the original nonlinear Riccati equation. Note that both assumptions admit explosion at time $T$. 
\end{remark}

Consider the following forward evolution equation: given the initial data $(t,x)$, 
\begin{equation} \label{forward evol}
dY_s=(A+A_0(s))Y_sds+\sum_{i=1}^\infty {\widehat C}_i(s)Y_sd\beta_s^i, \quad s\in (t,T].
\end{equation}

\begin{lemma} \label{LemmaAC} Let Assumptions~\eqref{A_0} and~\eqref{AC} hold true. There is a unique mild solution to~\eqref{forward evol} satisfying 
$$	\sup_{t\le s\le T} \mathbb{E}|Y_s|^2_{H}\le C |x|^2 $$
for a positive constant $C$.
\end{lemma}

\Dim 
First we prove the uniqueness.  Consider two solutions $Y^1$ and $Y^2$. Define $\Delta Y:=Y^1-Y^2.$
We have
$$
\Delta Y_s=\int_t^se^{A(s-r)}A_0(r)\Delta Y_r dr+\sum_{i=1}^\infty \int_t^se^{A(s-r)}{\widehat C}_i(r)\Delta Y_r d\beta_r^i, \quad  s\in [t,T]
$$
and 
\begin{eqnarray}
\mathbb{E}[|\Delta Y_s|^2]&\le& 2 \mathbb{E}\left|\int_t^se^{A(s-r)}A_0(r)\Delta Y_r dr\right|^2+2\mathbb{E}\left|\sum_{i=1}^\infty \int_t^se^{A(s-r)}{\widehat C}_i(r)\Delta Y_r d\beta_r^i\right|^2\nonumber \\
&\le &C\int_t^s \left((s-r)^{-2\alpha}+(T-r)^{-2\alpha}\right)\mathbb{E}[|\Delta Y_r|^2]\, dr\nonumber  \\
&\le& 2C\int_t^s(s-r)^{-2\alpha}\mathbb{E}[|\Delta Y_r|^2]\, dr. 
\end{eqnarray}
Thus, $\mathbb{E}[|\Delta Y_r|^2]=0$, and the uniqueness is proved. 

Then we prove the existence. 
Define by Picard's iteration: $Y^0\equiv 0$, and for $n\ge 0$, 
\begin{eqnarray*}
Y_s^{n+1}&=&e^{A(s-t)}x+\int_t^se^{A(s-r)}A_0(r)Y_r^n\, dr+\sum_{i=1}^\infty\int_t^se^{A(s-r)}{\widehat C}_i(r)Y_r^n\, d\beta_r^i.
\end{eqnarray*}
Thus, we have
\begin{eqnarray*}
\mathbb{E}\left[|Y_s^{n+1}|^2\right]&\le& 3|e^{A(s-t)}x|^2+2C\int_t^s(s-r)^{-2\alpha} E\left[|Y_r^n|^2\right]\, dr. \nonumber
\end{eqnarray*}
Denote by $\gamma$ the solution of the following integral equation:
\begin{equation}
\gamma_s=3|e^{A(s-t)}x|^2+2C\int_t^s (s-r)^{-2\alpha}\gamma_r\, dr, \quad s\in (t,T].
\end{equation}
By recurrence, we have $\mathbb{E}\left[|Y_s^n|^2\right]\le \gamma_s,$ for $s\in [t,T]$. 

Now we show that $\{Y^n, n\ge 0\}$ is a Cauchy sequence in $C_\calf([t,T];L^2(\Omega, H))$. We have 
\begin{eqnarray}
Y_s^{n+k+1}-Y_s^{n+1}&=&\int_t^se^{A(s-r)}A_0(r)\left(Y_r^{n+k}-Y_r^n\right)\, dr\nonumber \\
&&+\sum_{i=1}^\infty\int_t^se^{A(s-r)}{\widehat C}_i(r)\left(Y_r^{n+k}-Y_r^n\right)\, d\beta_r^i,\nonumber \\
\mathbb{E}\left[|Y_s^{n+k+1}-Y_s^{n+1}|^2\right]&\le&2C\int_t^s(s-r)^{-2\alpha} \mathbb E\left[|Y_r^{n+k}-Y_r^n|^2\right]\, dr. \nonumber
\end{eqnarray}

Define
$$
\phi_s=\limsup_{n}\sup_{k}\sup_{t\le r\le s}\mathbb{E}\left[|Y_r^{n+k+1}-Y_r^{n+1}|^2\right]. 
$$
We have 
\begin{eqnarray}
\sup_{t\le r \le s} \mathbb{E}\left[|Y_r^{n+k+1}-Y_r^{n+1}|^2\right]&\le&2C\int_t^s(s-r)^{-2\alpha}  \mathbb{E}\left[|Y_r^{n+k}-Y_r^n|^2\right]\, dr, \\
\phi_s&\le&2C\int_t^s(s-r)^{-2\alpha} \phi_r\, dr. 
\end{eqnarray}
This shows that $\phi=0$ and $\{Y^n\}$ is a Cauchy sequence in $C_\calf([0,T]; L^2(\Omega, H))$, and the existence of solution is proved. 
\finedim

\subsection{Lyapunov equation}

Let $G\in \cals(H)$ and $f\in L^1(0,T; \cals(H))$.  Assume that for  $\alpha\in (0,\frac{1}{ 2})$, 
\begin{equation}\label{f}
|f(s)|_{\call{(H)}}\le c(T-s)^{-2\alpha}. 
\end{equation}
Consider the following form of Lyapunov equation
\begin{equation}\label{lyapunov}
\left\{\begin{array}{l}
\displaystyle P'_t+ A^*P_t+P_tA +A_0^*(t)P_t+P_tA_0(t)+ \sum_{i=1}^\infty{\widehat C}_i^*(t)P_t{\widehat C}_i(t)+f_t=0, \quad t\in [0,T);   \\
P_T=G.
\end{array}\right.
\end{equation}

We look for a mild solution:
\begin{eqnarray}\label{LYA}
P_t&=&e^{A^*(T-t)}Ge^{A(T-t)}+\int_t^T e^{A^*(s-t)}f_se^{A(s-t)}ds \\
&&+\int_t^T e^{A^*(s-t)}\left[A_0^*(s)P_s+P_sA_0(s)+\sum_{i=1}^\infty {\widehat C}_i^*(s)P_s{\widehat C}_i(s)\right]e^{A(s-t)}ds. \nonumber
\end{eqnarray}


Using Yosida's approximation, we can prove that the following Lyapunov equation  (associated to a finite number of Brownian motions)
\begin{eqnarray}\label{LYA2}
P^n_t&=&e^{A^*(T-t)}Ge^{A(T-t)}+\int_t^T e^{A^*(s-t)}f_se^{A(s-t)}ds \\
&& +\int_t^T e^{A^*(s-t)}\left[A_0^*(s)P_s^n+P_s^nA_0(s)+\sum_{i=1}^n {\widehat C}_i^*(s)P^n_s{\widehat C}_i(s)\right]e^{A(s-t)}ds, \nonumber
\end{eqnarray}
has a unique solution $P^n\in C_s([0,T],\cals(H))$ (see, e.g. Da Prato~\cite{DaPrato1984}).

\begin{proposition}  \label{Prop3.3} Let Assumptions~\eqref{A_0}, \eqref{AC} and \eqref{f} hold true.
	Then, $P^n$  converges weakly to a bounded solution $P\in C_s([0,T];\cals(H))$ of \eqref{LYA} satisfying the estimate for some positive constant $C$, 
$$\sum_{i=1}^\infty \left|\langle {\widehat C}_i(s)^*P_s{\widehat C}_i(s)x, x\rangle \right|\le C (T-s)^{-2\alpha}|x|^2, \quad  s \in [0,T).$$
Moreover, we have the following representation of $P$:
\begin{equation}\label{eq:P}
\langle P_tx,x\rangle=\mathbb E\left[\langle GY_T^{t,x},Y_T^{t,x}\rangle+\int_t^T \langle f_sY_s^{t,x},Y_s^{t,x}\rangle ds\right], \quad t\in [0,T]
\end{equation}
where $Y^{t,x}$ is the mild solution  to  \eqref{forward evol}.
\end{proposition}

\Dim 
For each interger $n$, let $Y^{n, t,x}$ be the mild solution of 
$$dY_s=\left[A+A_0(s)\right]Y_sds+\sum_{i=1}^n {\widehat C}_i(s)Y_sd\beta_s^i, \quad s\in (t,T]; \quad Y_t=x.$$
We have the following representation:
\begin{equation}\label{eq:Pn}
\langle P^n_tx,x\rangle=\mathbb E\left[\langle GY^{n,t,x}_T,Y^{n, t, x}_T\rangle+\int_t^T \langle f_sY^{n,t,x}_s,Y^{n,t,x}_s\rangle ds\right]. 
\end{equation}
Since 
$$\lim_{n\to +\infty} \sup_{t\le s\le T} \mathbb E[|Y_s^{n,t,x}-Y_s^{t,x}|^2]= 0 $$ where $Y^{t,x}$ is the mild solution  to  \eqref{forward evol}, there exists $P_t\in \cals(H)$ such that $P_t^n$ converges to $P_t$ weakly and  we have by passing to the limit in \eqref{eq:Pn}  the desired representation~\eqref{eq:P}.

Set $z_i={\widehat C}_i(t)x.$

Let us estimate $\sum_{i=1}^\infty\left|\langle P_tz_i,z_i\rangle\right|$.
We have
\begin{eqnarray}
	\left|\langle P_tz_i,z_i\rangle\right|&=&\left|\mathbb E\left[\langle GY^{t,z_i}_T,Y^{t,z_i}_T\rangle+\int_t^T \langle f_s Y^{t,z_i}_s,Y^{t,z_i}_s\rangle ds\right]\right|\nonumber\\
	&\le& C\mathbb E[||Y^{t,z_i}_T||^2]+C\int_t^T (T-s)^{-2\alpha}\mathbb E[||Y^{t,z_i}_s||^2]ds.\label{Pzz}
\end{eqnarray}

As $Y^{t,z_i}$ is the mild solution of the following equation
$$Y_s^{t,z_i}=e^{A(s-t)}z_i+\int_t^s e^{A(s-r)} A_0(r)Y^{t,z_i}_rdr+\int_t^s e^{A(s-r)}\sum_{i=1}^\infty {\widehat C}_i(r)Y^{t,z_i}_rd\beta_r^i, $$
we have
\begin{eqnarray}\label{zi}
\mathbb E[||Y_s^{t,z_i}||^2]&\le& C||e^{A(s-t)}z_i||^2+C \int_t^s \left((s-r)^{-2\alpha}+(T-r)^{-2\alpha}\right)\mathbb E[||Y_r^{t,z_i}||^2]dr\nonumber\\
 &\le& C||e^{A(s-t)}z_i||^2+2C \int_t^s (s-r)^{-2\alpha}\mathbb E[||Y_r^{t,z_i}||^2]dr.
\end{eqnarray}

Note that
$$ \sum_{i=1}^n ||e^{A(s-t)}z_i||^2=\sum_{i=1}^n ||e^{A(s-t)}{\widehat C}_i(t)x||^2\le 
C\left((s-t)^{-2\alpha}+(T-t)^{-2\alpha}\right)||x||^2\le 2C(s-t)^{-2\alpha}||x||^2.$$

Finally, we get from~\eqref{zi} that 
$$\sum_{i=1}^n\mathbb E[||Y_s^{t,z_i}||^2]\le C(s-t)^{-2\alpha}||x||^2+C \int_t^s (s-r)^{-2\alpha}\sum_{i=1}^n\mathbb E[||Y_r^{t,z_i}||^2]dr.$$

By the generalized Gronwall's inequality (see Henry~\cite{Henry}), we have 

\begin{eqnarray*}
	\sum_{i=1}^n\mathbb E[||Y_s^{t,z_i}||^2]&\le& C(s-t)^{-2\alpha}||x||^2,
\end{eqnarray*}
and then letting $n\to \infty$, we have 
\begin{equation}\label{Ytz}
\sum_{i=1}^\infty\mathbb E[||Y_s^{t,z_i}||^2]\le C(s-t)^{-2\alpha}||x||^2.
\end{equation}

Furthermore from~\eqref{Pzz}  and ~\eqref{Ytz}, we have
\begin{eqnarray*}
\sum_{i=1}^\infty \left|\langle {\widehat C}_i(t)^*P_t{\widehat C}_i(t)x, x\rangle\right| &=&\sum_{i=1}^\infty \left|\langle P_tz_i,z_i\rangle\right| \\
&\le& C(T-t)^{-2\alpha}+C\int_t^T(T-s)^{-2\alpha}(s-t)^{-2\alpha}|x|_H^2\, ds\\
&=& C(T-t)^{-2\alpha}+C\int_0^1 (T-t)^{-2\alpha}(1-r)^{-2\alpha}(T-t)^{-2\alpha}r^{-2\alpha}(T-t)dr\\
&=&C(T-t)^{-2\alpha}+C(T-t)^{1-4\alpha}\le C(1+T^{1-2\alpha}) (T-t)^{-2\alpha}. 
\end{eqnarray*}

Passing to the limit in \eqref{LYA2} by letting $n\to \infty$, we prove that $P$ is the solution to~\eqref{LYA}.
\finedim

\begin{theorem} \label{Th3.4}
	There exists a unique solution $P\in C_s([0,T],\cals(H))$ for \eqref{lyapunov} such that 
	$$\sum_{i=1}^\infty \left|\langle{\widehat C}_i(s)^*P_s{\widehat C}_i(s)x, x\rangle \right|\le C (T-s)^{-2\alpha}|x|^2.$$
\end{theorem}

\Dim The existence of solution is already proved in the preceding proposition. Now we prove the uniqueness. 

Let $\tilde{P}$ be a solution, then it satisfies the following truncated Riccati equation:
\begin{eqnarray}
\tilde{P}_t&=&e^{A^*(T-t)}Ge^{A(T-t)}+\int_t^T e^{A^*(s-t)}\sum_{i=1}^n{\widehat  C}_i^*(s)\tilde{P}_s{\widehat C}_i(s)e^{A(s-t)}\, ds\nonumber\\
&&+\int_t^T e^{A^*(s-t)}(f_s+\sum_{i=n+1}^\infty {\widehat C}_i^*(s)\tilde{P}_s{\widehat C}_i(s))e^{A(s-t)}\, ds.
\end{eqnarray}

We have from~\eqref{eq:Pn}  the following representation 
$$\langle \tilde{P}_tx,x\rangle=\mathbb E\left[\langle GY^{n,t,x}_T,Y^{n,t,x}_T\rangle+\int_t^T \langle (f_s+\sum_{i=n+1}^\infty {\widehat C}_i^*(s)\tilde{P}_s{\widehat C}_i(s))Y^{n,t,x}_s,Y^{n,t,x}_s\rangle ds\right].$$
By passing to the limit, we deduce
$$\langle \tilde{P}_tx,x\rangle=\mathbb E\left[\langle GY_T^{t,x},Y_T^{t,x}\rangle+\int_t^T \langle f_s Y_s^{t,x},Y_s^{t,x}\rangle ds\right],$$
from which we deduce the uniqueness. \finedim

From~\eqref{eq:P}, we deduce also the following a priori estimate. 

\begin{proposition} \label{a priori} Let $P\in C_s([0,T],\cals(H))$ be the unique solution,
	then the following a priori estimate holds:
	$$|P_t|\le C (|G|_{\call{(H)}}+\int_t^T |f_s|_{\call{(H)}}ds).$$
\end{proposition}

\section{Riccati equation: existence and uniqueness of solutions}

In this section, we study the  Riccati equation associated to the linear-quadratic optimal control problem \eqref{eq:state mild} and \eqref{eq:cost 1}. Let us first state a lemma which will be used later. 

\begin{lemma} Let assumption~\eqref{control} hold true. For  $P\in \cals^+(H)$ such that for any $x\in H$, $$\sum_{j=1}^\infty\langle C_j^*(t)PC_j(t)x, x\rangle<\infty.$$
	Then, $\sum_{j=1}^N D_j^*(t)PC_j(t)$ converges strongly, whose limit is denoted by $\sum_{j=1}^\infty D_j^*(t)PC_j(t)$ and satisfies the  following estimate:
	$$
	\left|\sum_{j=1}^\infty D_j^*(t)PC_j(t)x\right|_U\le C \left(\sum_{j=1}^\infty \langle C_j^*(t)PC_j(t)x, x\rangle \right)^{\frac{1}{2}}, \quad x\in H
	$$
	for some constant $C>0$. 
\end{lemma}

\Dim
In view of Assumption~\eqref{control}, 
\begin{eqnarray*}
	\left|\sum_{j=M+1}^N D_j^*(t)PC_j(t)x\right|_U&=& \sup_{|y|_U\le 1} \left\langle\sum_{j=M+1}^ND_j^*(t)PC_j(t)x, y\right \rangle \\
	&\le &\sup_{|y|_U\le 1} \left(\sum_{j=M+1}^N |P^{\frac{1}{2}}C_j(t)x|_H^2\right)^{\frac{1}{2}}\left(\sum_{j=M+1}^N|P^{\frac{1}{2}}D_j(t)y|_H^2\right)^{\frac{1}{2}}\\
	&\le& C \left(\sum_{j=M+1}^N |P^{\frac{1}{2}}C_j(t)x|_H^2\right)^{\frac{1}{2}}\\
	&=&C \left(\sum_{j=M+1}^N \langle C_j^*(t)PC_j(t)x, x\rangle \right)^{\frac{1}{2}}.
\end{eqnarray*}
Hence the sequence $\sum_{j=1}^N D_j^*(t)PC_j(t)x$ is a Cauchy one, and we have the desired result. 
\finedim

Define for  $P\in \cals^+(H)$, 
$$
\Lambda(t, P):= R_t+\sum_{i=1}^\infty D_i^*(t)PD_i(t).
$$
Since $\Lambda(t, P)\ge \delta I_U$, we see that $\Lambda(t, P)$ has an inverse $\Lambda(t, P)^{-1}\le \frac{1}{\delta} I_U.$


Define for $s\in [0,T]$ and $P\in \cals^+(H)$, 
$$\lambda(s, P):=-\Lambda(s, P)^{-1}\left(B^*_sP+\sum_iD_i^*(s)PC_i(s)\right),$$
and  for $P\in \cals^+(H)$ such that $\sum_{j=1}^\infty \langle C_j^*(t)PC_j(t)x,x\rangle <\infty$ for each $x\in H$, 
\begin{eqnarray*}
	\hat{B}(t, P)&:=&-B_t
	\Lambda(t, P)^{-1}\left(B^*_tP+\sum_jD_j^*(t)PC_j(t)\right)=B_t\lambda(t, P),\\
	\hat{C}_i(t, P)&:=&C_i(t)-D_i(t)
	\Lambda(t, P)^{-1}\left(B^*_tP+\sum_jD_j^*(t)PC_j(t)\right)=C_i(t)+D_i(t)\lambda(t, P). 
\end{eqnarray*}

We have 

\begin{lemma}\label{Lipsch}  For $P\in C_s([0,T],\cals^+(H))$ such that 
	$$
	|\sum_{i=1} ^\infty C_i^*(s)P_s C_i(s)|_{\call(H)}\le C (T-s)^{-2\alpha}, 
	$$
we have
	$$	\sum_{i} |e^{At} {\widehat C}_i(s,P_s)x|^2\le c\left(t^{-2\alpha}+(T-s)^{-2\alpha}\right)|x|^2,$$
	$$\left|Q_s+\lambda^*(s, P_s)R_s\lambda(s, P_s)\right|_{\call{(H)}}\le c(T-s)^{-2\alpha},$$
	$$|{\widehat B} (s, P_s)|_{\call{(H)}}\le c(T-s)^{-\alpha}. $$
\end{lemma}

\Dim The third inequality is obvious. We now prove the first inequality. 
\begin{eqnarray}
&&	\sum_{i} |e^{At} {\widehat C}_i(s,P_s)x|^2\nonumber \\
&\le& 2 	\sum_{i} |e^{At} C_i(s)x|^2 +2\sum_{i} \left |e^{At} D_i(s)\Lambda (s, P_s)^{-1} \left(B^*_sP_s+\sum_jD_j^*(s)P_sC_j(s)\right)x\right |^2\nonumber \\
&\le& 2ct^{-2\alpha} |x|^2+2c|e^{At}|_{\call{(H)}}\left|\Lambda (s, P_s)^{-1} \left(B^*_sP_s+\sum_jD_j^*(s)P_sC_j(s)\right)x\right |^2\nonumber \\
&\le& 2ct^{-2\alpha} |x|^2+C|x|^2+C\left\langle \sum_j C_j^*(s)P_sC_j(s)x, x\right\rangle \nonumber \\
&\le& c\left(t^{-2\alpha} +(T-s)^{-2\alpha}\right)|x|^2.
\end{eqnarray}

It remains to  prove the second inequality. 
We have for each $x\in H$, since $R_s\le \Lambda(s, P_s)$, 
\begin{eqnarray*}
	&&\left\langle \lambda^*(s, P_s)R_s\lambda(s, P_s)x, x\right\rangle\\
	& \le&\left\langle (P_sB_s+\sum_iC_i^*(s)P_sD_i(s))\Lambda(s, P_s)^{-1}(B^*_sP_s+\sum_iD_i^*(s)P_sC_i(s))x, x\right\rangle\\
	&\le&
	2\left\langle P_sB_s\Lambda(s, P_s)^{-1}B^*_sP_sx, x\right\rangle +2\left\langle \sum_jC_j^*(s)P_sD_j(s)\Lambda(s, P_s)^{-1}\sum_iD_i^*(s)P_sC_i (s)x, x\right\rangle\\
	&\le& 2\left\langle P_sB_s\Lambda(s, P_s)^{-1}B^*_sP_sx, x\right\rangle  +2\left\langle \Lambda(s, P_s)^{-1}\sum_iD_i^*(s)P_sC_i (s)x, \sum_iD_i^*(s)P_sC_i (s)x\right\rangle\\
	&\le & c(T-s)^{-2\alpha}|x|_H^2. 
\end{eqnarray*}
\finedim

Let us consider the general Riccati equation:
\begin{eqnarray}
	P_t&=&e^{A^*(T-t)}Ge^{A(T-t)}+\int_t^T e^{A^*(s-t)}\sum_{i=1}^\infty C_i^*(s)P_sC_i(s)e^{A(s-t)}ds\nonumber\\
	& &+\int_t^T e^{A^*(s-t)}\left(Q_s-\lambda^*(s, P_s)\Lambda(s, P_s)\lambda(s, P_s)\right)e^{A(s-t)}ds, \quad t\in [0,T].\label{Riccati 1}
\end{eqnarray}
It is equivalent to the following form:
\begin{eqnarray}
P_t&=&e^{A^*(T-t)}Ge^{A(T-t)}+\int_t^T e^{A^*(s-t)} \left(\hat{B}^*(s, P_s)P_s+
P_s\hat{B}(s, P_s)\right)e^{A(s-t)}ds\nonumber\\
& &+\int_t^T e^{A^*(s-t)}\sum_{i=1}^\infty \hat{C}_i^*(s, P_s)P_s\hat{C}_i(s, P_s)e^{A(s-t)}ds \nonumber\\
& &+\int_t^T e^{A^*(s-t)}\left(Q_s+\lambda^*(s, P_s)R_s\lambda(s, P_s)\right)e^{A(s-t)}\, ds, \quad t\in [0,T]. \label{Riccati 2}
\end{eqnarray}

Our existence  proof will make  use of the following quasi-linearized sequence $\{P^N\}$ defined by the following Lyapunov equations: $P^0\equiv 0$,   and 
\begin{eqnarray}\label{PN+1}
	P^{N+1}_t&=&e^{A^*(T-t)}Ge^{A(T-t)}+\int_t^T e^{A^*(s-t)} \left(\hat{B}^*(s, P_s^N)P^{N+1}_s+
	P^{N+1}_s\hat{B}(s, P_s^N)\right)e^{A(s-t)}ds\nonumber\\
	& &+\int_t^T e^{A^*(s-t)}\sum_{i=1}^\infty \hat{C}_i^*(s, P_s^N)P^{N+1}_s\hat{C}_i(s, P_s^N)e^{A(s-t)}ds\nonumber \\
	& &+\int_t^T e^{A^*(s-t)}\left(Q_s+\lambda^*(s, P^N_s)R_s\lambda(s, P^N_s)\right)e^{A(s-t)}\, ds, \quad N=0,1,\ldots.
\end{eqnarray}

Note that if $P^N\in C_s([0,T],\cals^+(H))$ and satisfies the inequality for a positive constant $c$ which might depend on $N$: 
$$
|\sum_{i}C_i^*P^N_sC_i(s)|_{\call(H)}\le c (T-s)^{-2\alpha}, 
$$
then we see  from Lemma~\ref{Lipsch} and Theorem~\ref{Th3.4} that the preceding Lyapunov equation~\eqref{PN+1} has a unique solution $P^{N+1}\in C_s([0,T], \cals^+(H))$  satisfying also the last inequality. Since obviously $P^0$ satisfies the last inequality,  we can define  by induction a  sequence $P^{N+1}$  satisfying Lyapunov equation~\eqref{PN+1} for $N\ge 0$ .


\begin{lemma} The sequence 
	$\{P^N_t, N\ge 1\}$ is a non-increasing sequence of  self-adjoint operators for each $t\in [0,T]$.
\end{lemma}

\Dim
 Now we show that $P^N_t\ge P^{N+1}_t$ for $N\ge 1$.  

Define 
$\Delta P^N_t:=P^N_t-P^{N+1}_t, \quad  t\in [0,T]. $ We have 
\begin{eqnarray}\label{DeltaPN}
	\Delta P^N_t&=&\int_t^T e^{A^*(s-t)} \left(\hat{B}^*(s, P_s^{N-1})P^N_s+
	P^N_s\hat{B}(s, P_s^{N-1})\right)e^{A(s-t)}ds\nonumber\\
	& &+\int_t^T e^{A^*(s-t)}\sum_{i=1}^\infty \hat{C}_i^*(s, P_s^{N-1})P^N_s\hat{C}_i(s, P_s^{N-1})e^{A(s-t)}ds \nonumber\\
	& &+\int_t^T e^{A^*(s-t)}\lambda^*(s, P^{N-1}_s)R_s\lambda(s, P^{N-1}_s)e^{A(s-t)}\, ds\nonumber\\
	&&-\int_t^T e^{A^*(s-t)} \left(\hat{B}^*(s, P_s^N)P^{N+1}_s+
	P^{N+1}_s\hat{B}(s, P_s^N)\right)e^{A(s-t)}ds\nonumber\\
	& &-\int_t^T e^{A^*(s-t)}\sum_{i=1}^\infty \hat{C}_i^*(s, P_s^N)P^{N+1}_s\hat{C}_i(s, P_s^N)e^{A(s-t)}ds \nonumber\\
	& &-\int_t^T e^{A^*(s-t)}\lambda^*(s, P^N_s)R_s\lambda(s, P^N_s)e^{A(s-t)}\, ds. 
\end{eqnarray}
Define for $K\in \call (H, U)$ and  $P\in \cals^+(H)$ such that  $\sum_{i=1}^\infty\left \langle C_i^*(s)PC_i(s)x, x\right\rangle <\infty$, 
\begin{eqnarray}
F(s,K,P):=(B_sK)^*P+PB_sK+\sum_{i=1}^\infty\left [C_i(s)+D_i(s)K\right]^*P\left[C_i(s)+D_i(s)K\right]+K^*R_sK. 
\end{eqnarray}
We have for $K\in L(H, U)$,
$$
F(s,K,P)= F(s, \lambda(s, P), P) +\left[K-\lambda(s,P)\right]^*\Lambda(s,P)\left[K-\lambda(s,P)\right]\ge F(s, \lambda(s, P), P). 
$$

Equality~\eqref{DeltaPN} can be written into the following form:
\begin{eqnarray}
	\Delta P^N_t&=&\int_t^T e^{A^*(s-t)} \left(\hat{B}^*(s, P_s^N)	\Delta P^N_s+
	\Delta P^N_s\hat{B}(s, P_s^N)\right)e^{A(s-t)}ds\nonumber\\
	& &+\int_t^T e^{A^*(s-t)}\sum_{i=1}^\infty \hat{C}_i(s, P_s^N)	\Delta P^N_s\hat{C}_i(s, P_s^N)e^{A(s-t)}ds\nonumber\\
	&&+\int_t^T e^{A^*(s-t)} \left[F(s, \lambda(s, P^{N-1}_s), P_s^N) -F(s, \lambda(s, P^N_s), P_s^N)\right]e^{A(s-t)}ds. 
\end{eqnarray}
Note that $F(s, \lambda(s, P^{N-1}_s), P_s^N) -F(s, \lambda(s, P^N_s), P_s^N)\in \cals^+(H)$ for each $s\in [0,T]$. Therefore, we have from the representation theorem that $\Delta P^N_s\ge 0$. 
\finedim

We have the following theorem.  

\begin{theorem}\label{Th4.4} The Riccati equation~\eqref{Riccati 1} has a unique  solution $P\in C_s([0,T];\cals^+(H))$ such that 
	$$
	|\sum_{i}C_i^*P_sC_i(s)|_{\call(H)}\le C (T-s)^{-2\alpha}, \quad s\in [0,T]. 
	$$
\end{theorem}

{\bf Proof. }  First,  we see  from the last lemma that $P^N$ is a nondecreasing sequence of self-adjoint operators. Moreover, we see from~\eqref{eq:P} that each $P^N$ is non-negative. Using the monotone sequence theorem (see Kantorovich and Akilov~\cite[Theorem 1, p.  169]{KantorovichAkilov}), we see that $P^N_t$ converges strongly to a non-negative self-adjoint operator,  denoted by $P_t$, which also satisfies the last inequality.

As $P^N_t$ converges strongly to $P_t$, noting the following 
\begin{eqnarray*}
\Lambda(t, P^N_t)^{-1}-\Lambda(t, P_t)^{-1}&=&\Lambda(t, P^N_t)^{-1}\left(\Lambda(t, P_t)-\Lambda(t, P^N_t)\right)\Lambda(t, P_t)^{-1}\\
&=& \Lambda(t, P^N_t)^{-1}\left(\sum_{i=1}^\infty D_i^*(t)(P_t-P^N_t)D_i(t)\right)\Lambda(t, P_t)^{-1}, 
\end{eqnarray*}
we see that  $\Lambda(t, P^N_t)^{-1}$ converges strongly to $\Lambda(t, P_t)^{-1}$.

In view of our assumption~\eqref{control}, we have
\begin{eqnarray*}
&&\left|\sum_{j=1}^\infty D_j^*(t)P^N_tC_j(t)x-\sum_{j=1}^\infty D_j^*(t)P_tC_j(t) x\right|_U\nonumber \\
&=& \sup_{|y|_U\le 1} \left\langle \sum_{j=1}^\infty D_j^*(t)P^N_tC_j(t)x-\sum_{j=1}^\infty D_j^*(t)P_tC_j(t)x, y\right \rangle \\
	&=& \sup_{|y|_U\le 1} \sum_{j=1}^\infty \left\langle (P^N_t-P_t)^{\frac{1}{2}}C_j(t)x, (P^N_t-P_t)^{\frac{1}{2}}D_j(t)y\right \rangle\\
	&\le & C \left(\sum_{j=1}^\infty\left |(P^N_t-P_t)^{\frac{1}{2}}C_j(t)x\right|^2\right)^{\frac{1}{2}}\\
	&= & C \left(\sum_{j=1}^\infty\left \langle C_j^*(t)(P^N_t-P_t)C_j(t)x, x\right\rangle\right)^{\frac{1}{2}}. 
\end{eqnarray*}
Since $\langle C_j^*(t)(P^N_t-P_t)C_j(t)x, x\rangle\le \langle C_j^*(t)(P^1_t-P_t)C_j(t)x, x\rangle$ and 
$$
\sum_{j=1}^\infty \langle (C_j^*(t)(P^1_t-P_t)C_j(t)x, x\rangle<\infty, 
$$
by the Dominated Convergence Theorem, we see that $\sum_{j=1}^\infty D_j^*(t)P^N_tC_j(t)$ converges strongly
to $\sum_{j=1}^\infty D_j^*(t)P_tC_j(t)$.  Therefore, the non-homogeneous term in the Lyapunov equation of $P^{N+1}$ converges strongly. 

By passing to the strong limit in the Lyapunov equation~\eqref{PN+1}, we conclude that $P$ is a solution. 

Finally, we show the uniqueness. Let $\widetilde{P}$ be another  solution of Riccati equation~\eqref{Riccati 1} such that 
$$
|\sum_{i}C_i^*(s){\widetilde P}_sC_i(s)|_{\call(H)}\le C (T-s)^{-2\alpha}, \quad s\in [0,T). 
$$
Define $\delta P:=P-\widetilde{P}$. Then proceeding identically as in the last lemma, we have 
\begin{eqnarray}
\delta P_t&=&\int_t^T e^{A^*(s-t)} \left(\hat{B}^*(s, P_s)	\delta P_s+
\delta P_s\hat{B}(s, P_s)\right)e^{A(s-t)}ds\nonumber\\
& &+\int_t^T e^{A^*(s-t)}\sum_{i=1}^\infty \hat{C}_i(s, P_s)	\delta P_s\hat{C}_i(s, P_s)e^{A(s-t)}ds\nonumber\\
&&+\int_t^T e^{A^*(s-t)} \left[ F(s, \lambda(s, P_s), {\widetilde P}_s)-F(s, \lambda(s, {\widetilde P}_s), {\widetilde P}_s)\right]e^{A(s-t)}ds. 
\end{eqnarray}
Since $F(s, \lambda(s, P_s), P_s)-F(s, \lambda(s, {\widetilde P}_s), P_s)$ is non-negative,  we have  $\delta P\ge 0$. By symmetry, we also have $\delta P\le 0$. Hence, we have $\delta P\equiv 0$. 
\finedim

\section{Optimal feedback control}
In this section, we study the linear quadratic optimal control problem~\eqref{control system} and \eqref{eq:cost 1}. 

Note that It\^o's formula could not be applied to systems driven by a white noise. To overcome the difficulty, we truncate the white noise  by a finite number of Brownian motions. 

Define $X^N$ to be the unique solution of the following truncated state equation:
\begin{eqnarray}\label{eq:stateN2}
dX^N_t&=&(AX^N_t+B_tu_t)\,dt+\sum_{j=1}^N \left[C_j(t)X^N_t+D_j(t)u_t\right]\,d \beta^j_t,\\
X^N_0&=&x\in H.\nonumber
\end{eqnarray}

We denote by $P^N$ the solution of the following truncated Lyapunov equation:
\begin{eqnarray}\label{truncation P}
P_t^N&=&e^{A^*(T-t)}Ge^{A(T-t)}+\int_t^T e^{A^*(s-t)} \left(\hat{B}^*(s, P_s)P_s^N+
P_s^N\hat{B}(s, P_s)\right)e^{A(s-t)}ds\nonumber\\
& &+\int_t^T e^{A^*(s-t)}\sum_{i=1}^N \hat{C}_i^*(s, P_s)P_s^N\hat{C}_i(s, P_s)e^{A(s-t)}ds\\
& &+\int_t^T e^{A^*(s-t)}\left(Q_s+\lambda^*(s, P_s)R_s\lambda(s, P_s)\right)e^{A(s-t)}\, ds, \quad t\in [0,T], \nonumber
\end{eqnarray}
where $P$ is the unique solution of the Riccati equation~\eqref{Riccati 1}.

We have 

\begin{lemma}\label{P convergence}  For $t\in [0,T]$, 
	$P^N_t$ is non-decreasing, and is bounded from above, and strongly converges to $P_t$. 
\end{lemma}

\Dim 
In view of~\eqref{truncation P},  we see from the representation that $P^N_t$ is nonnegative for each $N\ge 1$. Furthermore, we have 
\begin{eqnarray*}
P^{N+1}_t-P^N_t&=&\int_t^T e^{A^*(s-t)} \left(\hat{B}^*(s, P_s)(P_s^{N+1}-P_s^N)+
(P_s^{N+1}-P_s^N)\hat{B}(s, P_s)\right)e^{A(s-t)}ds\nonumber\\
& &+\int_t^T e^{A^*(s-t)}\sum_{j=1}^N C_j^*(s)(P^{N+1}_s-P^N_s)C_j(s)e^{A(s-t)}ds\\
& &+\int_t^T e^{A^*(s-t)}{\widehat C}_{N+1}^*(s, P_s)P_s^{N+1}{\widehat C}_{N+1}(s, P_s)e^{A(s-t)}ds. 
\end{eqnarray*}

From the representation of the solution of Lyapunov equation, it is clear that 	$P^N$ is non-decreasing. 	Using the same argument to consider the equation of $P-P^N$,  we see that $P_t^N\le P_t$.  Hence there is a bounded ${\overline P}\le P$ such that $P^N$ strongly converges to ${\overline P}$ which satisfies the Lyapunov equation:
\begin{eqnarray}
{\overline P}_t&=&e^{A^*(T-t)}Ge^{A(T-t)}+\int_t^T e^{A^*(s-t)} \left(\hat{B}^*(s, P_s){\overline P}_s+
{\overline P}_s\hat{B}(s, P_s)\right)e^{A(s-t)}ds\nonumber\\
& &+\int_t^T e^{A^*(s-t)}\sum_{i=1}^\infty \hat{C}_i^*(s, P_s){\overline P}_s\hat{C}_i(s, P_s)e^{A(s-t)}ds\\
& &+\int_t^T e^{A^*(s-t)}\left(Q_s+\lambda^*(s, P_s)R_s\lambda(s, P_s)\right)e^{A(s-t)}\, ds, \quad t\in [0,T], \nonumber
\end{eqnarray}
and
$$
|\sum_{i}C_i^*{\overline P}_sC_i(s)|_{\call(H)}\le C (T-s)^{-2\alpha}, \quad s\in [0,T]. 
$$

Since $P$ (as a solution to the Riccati equation) is also a solution to the preceding Lyapunov equation with the non-homogeneous term being $f_s:= Q_s+\lambda^*(s, P_s)R_s\lambda(s, P_s)$,  we conclude from the uniqueness of the solution to the Lyapunov equation that ${\overline P}=P$. 
\finedim

\begin{theorem} \label{EnergyEquality}
	The cost functional has the following representation :
	$$J(x,u)=\langle P_0x,x\rangle +\mathbb E\left[\int_0^T \langle \Lambda (s, P_s) (u_s-\lambda(s, P_s)X_s), u_s-\lambda(s, P_s)X_s\rangle ds\right].$$
	 The following feedback form:
	\begin{equation}
	{\overline u}_t=\lambda (t,P_t) {\overline X}_t, \quad t\in [0,T],
	\end{equation}
	with ${\overline X}$ being the solution of the associated feedback system,  is admissible and optimal. 
	\end{theorem}

\Dim
We have the duality between the truncated state equation and the truncated Lyapunov equation by Yosida approximation of $A$:
\begin{eqnarray*}\label{dualityN2}
&&\mathbb E\left[\langle GX^N_T,X^N_T\rangle+\int_0^T \langle(Q_s+\lambda^*(s, P_s)R_s\lambda(s, P_s))X^N_s,X^N_s\rangle\right]ds\\
		&= &\langle P^N_0x,x\rangle+2\mathbb E\int_0^T\langle P^N_sX^N_s, B_su_s\rangle ds\\
	&&	-\mathbb E\int_0^T\left\langle \left[{\widehat B}^*(s, P_s)P^N_s+P_s^N{\widehat B}(s, P_s)+\sum_{j=1}^N{\widehat C}^*_j(s, P_s)P^N_s {\widehat C}_j(s, P_s)\right]X^N_s, X^N_s\right \rangle ds\\ 
	&& +\mathbb E\int_0^T\sum_{j=1}^N\left\langle P^N_s(C_j(s)X_s^N+D_j(s)u_s),  C_j(s)X_s^N+D_j(s)u_s\right\rangle  ds.
\end{eqnarray*}
Setting $N\to \infty$, noting the following limit
$$
\lim_{N\to \infty}\sup_{0\le t\le T}\mathbb{E}\left[|X^N_t-X_t|^2\right]=0,
$$
and Lemma~\ref{P convergence}, we have 
\begin{eqnarray*}\label{dualityN4}
	&&\mathbb E\left[\langle GX_T,X_T\rangle+\int_0^T \langle(Q_s+\lambda^*(s, P_s)R_s\lambda(s, P_s))X_s,X_s\rangle\right]ds\\
	&= &\langle P_0x,x\rangle+2\mathbb E\int_0^T\langle P_sX_s, B_su_s\rangle ds\\
	&&	-\mathbb E\int_0^T\left\langle \left[{\widehat B}^*(s, P_s)P_s+P_s{\widehat B}(s, P_s)+\sum_{j=1}^\infty{\widehat C}^*_j(s, P_s)P_s {\widehat C}_j(s, P_s)\right]X_s, X_s\right \rangle ds\\ 
	&& +\mathbb E\int_0^T\sum_{j=1}^\infty\left\langle P_s(C_j(s)X_s+D_j(s)u_s),  C_j(s)X_s+D_j(s)u_s\right\rangle  ds.
\end{eqnarray*}
Then, we have for any admissible control $u$, 
\begin{eqnarray*}\label{dualityN5}
	J(x,u)&=&\mathbb E\left[\langle GX_T,X_T\rangle+\int_0^T \langle Q_s X_s, X_s\rangle ds+\int_0^T \langle R_su_s, u_s\rangle ds\right]\\
	& =&\langle P_0 x, x\rangle+2\mathbb E\int_0^T\langle P_sX_s, B_su_s\rangle ds+\mathbb E\int_0^T \langle R_su_s, u_s\rangle ds\\
	&&-\int_0^T \langle\lambda^*(s, P_s)R_s\lambda(s, P_s)X_s,X_s\rangle ds\\
	&&	-\mathbb E\int_0^T\left\langle \left[{\widehat B}^*(s, P_s)P_s+P_s{\widehat B}(s, P_s)+\sum_{j=1}^\infty{\widehat C}^*_j(s, P_s)P_s {\widehat C}_j(s, P_s)\right]X_s, X_s\right \rangle ds\\ 
&& +\mathbb E\int_0^T\sum_{j=1}^\infty\left\langle P_s(C_j(s)X_s+D_j(s)u_s),  C_j(s)X_s+D_j(s)u_s\right\rangle  ds.
	\end{eqnarray*}
Noting that 
\begin{eqnarray*}
&&\sum_{j=1}^\infty C_j^*(s)PC_j(s)-\lambda^*(s, P)R_s\lambda(s, P)\\
&&-[{\widehat B}^*(s, P)P+P{\widehat B}(s, P)+\sum_{j=1}^\infty{\widehat C}^*_j(s, P)P {\widehat C}_j(s, P)]\\
&=&\lambda^*(s, P)\Lambda(s, P)\lambda(s, P)
\end{eqnarray*}
and
\begin{eqnarray*}
	2B_s^*P_s+2\sum_{j=1}^\infty D_j^*(s)P_sC_j(s)=-2\Lambda(s, P_s)\lambda(s, P_s), 
\end{eqnarray*}
we have for any admissible control $u$, 
\begin{eqnarray*}
J(x,u)	&=&\langle P_0x,x\rangle +\mathbb E\left[\int_0^T \langle \Lambda (s, P_s) (u_s-\lambda(s, P_s)X_s), u_s-\lambda(s, P_s)X_s\rangle ds\right]\\
&\ge& \langle P_0x,x\rangle .
\end{eqnarray*}

In view of Lemmas~\ref{Lipsch} and~\ref{LemmaAC},  we see that the closed-loop state equation has a unique solution ${\overline X}$, satisfying the following estimate
$$
\sup_{0\le t\le T}\mathbb{E}[|{\overline X_t}|^2_H]<\infty.
$$ 
Then
\begin{eqnarray*}
\mathbb{E}\int_0^T|{\overline u}_s|^2_U ds&=&\mathbb{E}\int_0^T|\lambda(s, P_s) {\overline X}_s|^2_U ds\\
&\le& C \int_0^T(T-s)^{-2\alpha}ds\sup_{0\le s\le T} \mathbb{E}|{\overline X}_s|^2<\infty.
\end{eqnarray*}
This shows that ${\overline u}$ is admissible and $J(x, {\overline u})=\langle P_0x,x\rangle $. 
Therefore, ${\overline u}$ is optimal. 
\finedim

\section{Algebraic Riccati equation} 

In this section, we discuss the solvability of algebraic Riccati equation.  For this, we need the following notion of stabilizability. Now we suppose that all the coefficients $B,C,D,Q,R$ are time-invariant. 

\begin{definition} We say that the system $(A,B,C,D)$ is feedback stabilizable if there is an operator $K\in \call(H, U)$ such that the system corresponding to the feedback control $u=KX$ is stable, i.e. for any initial state $x\in H$, 
	$$\mathbb{E} \int_0^\infty |X_t^{0,x}|^2dt<\infty. $$	
\end{definition}

\begin{theorem}\label{alg} Assume that the system $(A,B,C,D)$ is feedback stabilizable. Then there is a non-negative operator $P\in \cals^+(H)$ such that $\sum_{i=1}^\infty C_i^*PC_i\in \cals^+(H)$ and for any $T>0$, 
\begin{eqnarray}\label{Alg eq}
P&=&e^{A^*(T-t)}Pe^{A(T-t)}+\int_t^T e^{A^*(s-t)} \left(B^*P+
 PB+\sum_{i=1}^\infty C_i^*PC_i\right)e^{A(s-t)}ds\\
& &+\int_t^T e^{A^*(s-t)}\left[Q-\lambda^*(P)\left(R+\sum_{i=1}^\infty D_i^*PD_i\right)^{-1}\lambda(P)\right]e^{A(s-t)}\, ds, \quad t\in [0,T]. \nonumber
\end{eqnarray}
\end{theorem}

\Dim
Let $N\ge 1$. Consider the following Riccati equation on $[0,N]$:
\begin{eqnarray}
P^N_t&=&\int_t^N e^{A^*(s-t)} \left(B^*P^N_s+
P^N_sB+\sum_{i=1}^\infty C_i^*P^N_sC_i\right)e^{A(s-t)}ds\nonumber\\
& &+\int_t^N e^{A^*(s-t)}\left[Q-\lambda^*(P^N_s)\left(R+\sum_{i=1}^\infty D_i^*P^N_sD_i\right)^{-1}\lambda(P^N_s)\right]e^{A(s-t)}\, ds, \quad t\in [0,N]\nonumber
\end{eqnarray}
with the following estimate
$$
|\sum_{i=1}^\infty C_i^*P^N_sC_i|_{\call(H)}\le C_N(N-s)^{-2\alpha}. 
$$
It is easy to see that $P^N$ is non-decreasing in $N$. By the stabilizability assumption, there is a feedback control $u=KX$ such that $\mathbb{E}\int_0^\infty |X_s|^2\, ds< \infty$ and 
\begin{eqnarray}
\langle P^N_0 x, x\rangle &\le& \mathbb{E} \langle P^N_N X_N^{0,x}, X_N^{0,x}\rangle +\mathbb{E} \int_0^N \left(\langle QX_s^{0,x}, X_s^{0,x}\rangle +\langle Ru_s, u_s\rangle\right)\, ds\nonumber \\
&\le & \mathbb{E} \int_0^\infty \left(\langle QX_s^{0,x}, X_s^{0,x}\rangle +\langle Ru_s, u_s\rangle\right)\, ds=:C |x|^2 
\end{eqnarray}
with the number $C$ not depending on $N$.  Using the time-invariance of the underlying coefficients, we also have for each $t\in [0,\infty)$,  $\langle P_t^Nx, x\rangle \le C_t |x|^2$ with the number $C_t$ not depending on $N$.  Thus there exists $P_t$ such that $P^N_t$ converges to $P_t$ in a strong way.

For $t\le T\le N$, $P^N$ is the solution of the following Riccati equation 
\begin{eqnarray}\label{algP_n}
P^N_t&=&e^{A^*(T-t)}P^N_Te^{A(T-t)}+\int_t^T e^{A^*(s-t)} \left(B^*P^N_s+
P^N_sB+\sum_{i=1}^\infty C_i^*P^N_sC_i\right)e^{A(s-t)}ds\\
& &+\int_t^T e^{A^*(s-t)}\left[Q-\lambda^*(P^N_s)\left(R+\sum_{i=1}^\infty D_i^*P^N_sD_i\right)^{-1}\lambda(P^N_s)\right]e^{A(s-t)}\, ds, \quad t\in [0,T].\nonumber
\end{eqnarray}

As $|P^N_T|_{\call(H)}\le C_T$, from Theorem~\ref{Th3.4}, there exists a constant $C_T'$ such that
$$
|\sum_{i=1}^\infty C_i^*P^N_sC_i|_{\call(H)}\le C_T'(T-s)^{-2\alpha},\quad s\in [0,T),
$$
or equivalently 
$$
\sum_{i=1}^\infty \langle C_i^*P^N_sC_ix, x\rangle \le C_T'(T-s)^{-2\alpha}|x|^2,\quad (s,x)\in [0,T)\times H. 
$$
Passing to the limit in $N$, using Fatou's lemma, we derive
\begin{equation}\label{AC'}
|\sum_{i=1}^\infty C_i^*P_sC_i|_{\call(H)}\le C_T'(T-s)^{-2\alpha},\quad s\in [0,T).
\end{equation}

Taking the strong limit in \eqref{algP_n}, we deduce that $P_t$ satisfies the following equation:
\begin{eqnarray}\label{algP}
P_t&=&e^{A^*(T-t)}P_Te^{A(T-t)}+\int_t^T e^{A^*(s-t)} \left(B^*P_s+
P_sB+\sum_{i=1}^\infty C_i^*P_sC_i\right)e^{A(s-t)}ds\nonumber\\
& &+\int_t^T e^{A^*(s-t)}\left[Q-\lambda^*(P_s)\left(R+\sum_{i=1}^\infty D_i^*P_sD_i\right)^{-1}\lambda(P_s)\right]e^{A(s-t)}\, ds, \quad t\in [0,T]\nonumber
\end{eqnarray}

Due to the time invariance of the coefficients, we prove that $P_t$ does not depend on $t$.  From~\eqref{AC'}, we see that 
$\sum_{i=1}^\infty C_i^*PC_i\in \call(H)$, and that $P$ is the mild solution of~\eqref{alg}.
\finedim

\begin{theorem}\label{alg verification} Let $Q\in \cals^+(H)$ be positive. Assume that the non-negative operator $P\in \cals^+(H)$ satisfies  $\sum_{i=1}^\infty C_i^*PC_i\in \cals^+(H)$ and algebraic Riccati equation~\eqref{Alg eq}. 
	Then, the feedback law ${\overline u}= \lambda (P) {\overline X}$ is admissible and optimal, and the value function $=J(x,{\overline u})=\langle Px, x\rangle, x\in H.$ Consequently, the non-negative solution $P$ of algebraic Riccati equation~\eqref{Alg eq} such that $\sum_{i=1}^\infty C_i^*PC_i\in \cals^+(H)$  is unique. 
\end{theorem}

\Dim
For any admissible $u$, there is a sequence  $T_i\to \infty$ such that 
$$\lim_{i\to \infty}\mathbb{E}[|X_{T_i}|^2]=0.$$
Since $P_t\equiv P$ solves the Riccati equation on the finite time interval $[0,T_i]$ with the terminal condition $P$, we have 
$$
\mathbb{E} \langle PX_{T_i}, X_{T_i}\rangle +\mathbb{E} \int_0^{T_i}\left(\langle QX_s, X_s\rangle +\langle Ru_s, u_s\rangle \right)ds \ge \langle Px,x \rangle. 
$$
Letting $i\to \infty$, we have $J(x,u)\ge \langle Px,x \rangle$. 

Now we  prove that ${\overline u}$ is admissible and $J(x,{\overline u})= \langle Px,x \rangle$. Again since $P_t\equiv P$ solves the Riccati equation on the finite time interval $[0,T]$ with the terminal condition $P$, we have 
$$
\mathbb{E} \langle P{\overline X}_T, {\overline X}_T\rangle +\mathbb{E} \int_0^T\left(\langle Q{\overline X}_s, {\overline X}_s\rangle +\langle R{\overline u}_s, {\overline u}_s\rangle \right)ds =\langle Px,x \rangle. 
$$
By the monotone convergence theorem, we have 
$$
\mathbb{E} \int_0^\infty\left(\langle Q{\overline X}_s, {\overline X}_s\rangle +\langle R{\overline u}_s, {\overline u}_s\rangle \right)ds \le \langle Px,x \rangle. 
$$
As $Q$ is positive,  ${\overline u}$ is admissible and thus $J(x,{\overline u})= \langle Px,x \rangle. $
\finedim

\section{Null controllability of SPDEs via Riccati equation}

In this section, we characterize the null controllability of the system~\eqref{control system} via the existence of Riccati equation with the singular terminal value in the spirit of S\^{i}rbu and Tessitore~\cite{SirbuTessitore2001}. 

\begin{definition} The system~\eqref{control system} is $T$-null (exact) controllable if for any $(t, x)\in [0,T)\times H$,  there is $u\in L^2_{\calf}(t,T; U)$ such that $X^{t,x,u}_T=0$, $\mathbb{P}$-almost surely; and it is null controllable if it is  $T$-null controllable for each $T>0$. 
\end{definition}

For each $T>0$ and $x\in H$,  we consider the following optimal null-controllability control problem: 
\begin{equation}\label{cost S}
V(t,x):=\min_{u\in L^2_{\calf}(t,T; U)} J(t,x;u):=\mathbb{E}\int_t^T(|X_s^{t,x;u}|^2+|u_s|^2)\,ds
\end{equation}
subject to $X^{t,x;u}_T=0$. If there is no $u$ satisfying $X^{t,x;u}_T=0$, we set $V(t,x)=+\infty$.

Let $Id$ denote the identity operator in $H$. We introduce the following Riccati equation:
\begin{eqnarray}\label{Riccati S}
P_t&=&e^{A^*(T'-t)}P_{T'}e^{A(T'-t)}+\int_t^{T'} e^{A^*(s-t)}\sum_{i=1}^\infty C_i^*(s)P_sC_i(s)e^{A(s-t)}ds\nonumber\\
& &+\int_t^{T '}e^{A^*(s-t)}\left(Id-\lambda^*(s, P_s)\Lambda(s, P_s)\lambda(s, P_s)\right)e^{A(s-t)}ds, \quad 0\le t\le T'< T;
\end{eqnarray}
with the following singular terminal condition $P_T=+\infty$ in the following sense: for any $x\in H$ such that $x\not=0$,
\begin{equation}
\lim_{(t,y)\to (T,  x) }\langle P_ty, y\rangle =+\infty. 
\end{equation}

We then have 

\begin{theorem} \label{equivalence}
For given $T>0$,	the following conditions are equivalent: 
	
	(i) the Riccati equation~\eqref{Riccati S} has a mild solution $P$ satisfying the singular terminal condition $P_T=+\infty$ and the map $s\mapsto \sum_{i=1}^\infty C_i^*(s)P_sC_i(s)\in \call(H)$ is bounded in any compact interval of $[0, T)$;
	
	(ii) the state system~\eqref{control system} is $T$-null controllable. 
	
	If the system~\eqref{control system} is $T$-null controllable, then the associated optimal  null-control problem with the cost~\eqref{cost S} has the optimal control of the following feedback form:
	$$
	{\overline u}_s := \lambda(s, P_s){\overline X}_s, \quad s\in [0, T]
	$$
where $P$ is the solution of the Riccati equation~\eqref{Riccati S} with the singular terminal condition $P_T=+\infty$. 
\end{theorem}

\Dim First we prove that Assertion (i) implies (ii). In fact, if $P$ is the solution of the Riccati equation~\eqref{Riccati S} with the singular terminal condition $P_T=+\infty$, then its restriction on $[t, s]$  for $s\in (t, T)$ can be regarded as  the solution of the Riccati equation~\eqref{Riccati S} with the  terminal condition $P_s$.  Set $	{\overline u}_s := \lambda(s, P_s){\overline X}_s$ for $s\in [0, T]$. From Theorem~\ref{EnergyEquality}, we have 
\begin{eqnarray}
\langle P_t x, x\rangle&=&\mathbb{E} \langle P_s{\overline X}_s^{t,x;{\overline u}}, {\overline X}_s^{t,x;{\overline u}}\rangle+\mathbb{E}\int_t^s(|{\overline X}_r^{t,x;{\overline u}}|^2+|{\overline u}_r|^2)\,dr\nonumber \\
&\ge & \mathbb{E} \langle P_s{\overline X}_s^{t,x;{\overline u}}, {\overline X}_s^{t,x;{\overline u}}\rangle. \label{EngIne S}
\end{eqnarray}
Therefore, we have ${\overline u}\in L^2_\calf([t,T]; U)$ and we can extend ${\overline X}^{t,x;{\overline u}}$ to $[t, T]$ lying in $C_\calf([t,T];L^2(\Omega, H))$. From the inequality~\eqref{EngIne S} and Fatou's lemma, we have 
$$
\langle P_t x, x\rangle\ge  \mathbb{E} \liminf_{s\to T-} \left[\langle P_s{\overline X}_s^{t,x;{\overline u}}, {\overline X}_s^{t,x;{\overline u}}\rangle\right]\ge \mathbb{E}  \left[(+\infty) I_{{\overline X}_T^{t,x;{\overline u}}\not=0}\right] . 
$$
Hence, ${\overline X}^{t,x;{\overline u}}_T=0$, $\mathbb{P}$-almost surely. Assertion (ii) is proved. 

Now we show that Assertion (ii) implies (i). For any integer $n\ge 1$,  
\begin{equation}\label{cost Sn}
V^n(t,x):=\min_{u\in L^2_{\calf}(t,T; U)} J^n(t,x;u):=n\mathbb{E}|X_T^{t,x;u}|^2+ \mathbb{E}\int_t^T(|X_s^{t,x;u}|^2+|u_s|^2)\,ds. 
\end{equation}
It is associated to the following Riccati equation 
\begin{eqnarray}\label{Riccati Sn}
P_t&=&ne^{A^*(T-t)}e^{A(T-t)}+\int_t^{T} e^{A^*(s-t)}\sum_{i=1}^\infty C_i^*(s)P_sC_i(s)e^{A(s-t)}ds\nonumber\\
& &+\int_t^{T }e^{A^*(s-t)}\left(Id-\lambda^*(s, P_s)\Lambda(s, P_s)\lambda(s, P_s)\right)e^{A(s-t)}ds, \quad 0\le t\le T. 
\end{eqnarray}
Denoting by $P^n$ its unique solution, we have for $s<T$, 
\begin{eqnarray}\label{Riccati Sns}
P_t^n&=&e^{A^*(s-t)}P_s^ne^{A(s-t)}+\int_t^{s} e^{A^*(r-t)}\sum_{i=1}^\infty C_i^*(r)P_r^nC_i(r)e^{A(r-t)}dr\nonumber\\
& &+\int_t^{s }e^{A^*(r-t)}\left(Id-\lambda^*(r, P_r^n)\Lambda(r, P^n_r)\lambda(r, P_r^n)\right)e^{A(r-t)}dr, \quad 0\le t\le s. 
\end{eqnarray}

From Theorem~\ref{EnergyEquality}, we see that $P^n$ is non-decreasing in $n$. 
Moreover, from the $T$-null controllability, there is $u^0\in \call^2_\calf(t,T;U)$ such that $X_T^{t,x;u^0}=0.$ Hence, 
$$
\langle P^n_tx, x\rangle\le \mathbb{E}\int_t^T(|X_s^{t,x;u^0}|^2+|u_s^0|^2)\,ds.
$$
Consequently, the sequence $P^n_t$ has a strong limit in $\cals^+(H)$, which is denoted by $P_t$.  

For $s'\in (s,T)$, we have 
\begin{equation}
\left|\sum_{i=1}^\infty C_i^*(t)P_t^nC_i(t)\right|_{\call(H)}\le \frac{C_{s'}}{(s'-t)^{2\alpha}}\le \frac{C_{s'}}{(s'-s)^{2\alpha}}.
\end{equation}
Letting $n\to +\infty,$ we have 
$$
\left|\sum_{i=1}^\infty C_i^*(t)P_tC_i(t)\right|_{\call(H)}\le  \frac{C_{s'}}{(s'-s)^{2\alpha}}, 
$$
meaning that the sum is bounded in $\call(H)$. 

Taking the strong limit in~\eqref{Riccati Sns}, we see that $P$ is a mild solution of Riccati equation~\eqref{Riccati S} on the time interval $[0,T)$. 

Furthermore, we have for any integer $n$, 
$$
\liminf_{s\to T-,\  y\to x}\langle P_sy, y\rangle \ge \liminf_{s\to T-, \ y\to x}\langle P_s^ny, y\rangle =n|x|^2. 
$$
This shows that $P$ satisfies the singular terminal condition at time $T$. \finedim



\section{Examples: LQ optimal control of the Anderson model}

\begin{example} 
	Consider the following controlled Anderson model, that is, the following controlled stochastic heat equation in $[0,1]$: 
	\begin{eqnarray}
	dX_t(y)&=&\frac{\partial^2}{\partial^2 y}X_t(y)\, dt+b(t,y) u(t,y)\, dt+X_t(y) dW(t,y);\label{Anderson}\\
	X_t(0)&=&X_t(1)=0, \quad t\in [0,T];\\
	X_0(y)&=&x(y), \quad y\in [0,1]. 
	\end{eqnarray}
	The cost functional reads:
	\begin{eqnarray}
	J(x,u)=\mathbb E\int_0^T\int_0^1\left[q(t,y)X^2_t(y)+r(t,y)u^2(t,y)\right]\, dydt+\mathbb E\int_0^1g(y)X_T^2(y)\, dy. 
	\end{eqnarray}
	In the above example, $H=L^2(0,1)$, and  $W$ is an $H$-valued cylindrical Wiener process.  We choose an orthonormal basis $\{e_i, i=1,2,\ldots\}$ in the space $H$ such that
	$$
	\sup_i\sup_{y\in [0,1]} |e_i(y)|<\infty. 
	$$ 
	$A$ is the realization of the second derivative operator with the zero Dirichlet boundary conditions, and all the functions $b,q,r,g$ are measurable and  bounded . So $\cald(A)=H^2([0,1])\cap H_0^1([0,1])$ and $A\psi=\psi''$ for all $\psi\in \cald(A)$,  $C_i\phi(y):=e_i(y)\phi(y)$, and $(B_t\phi)(y)=b(t,y)\phi(y)$ for $\phi\in H$.  
	Then the pair $(A,C)$ satisfies (see Da Prato and Zabczyk~\cite{DaZa2}) the inequality~\eqref{AC0}. 
	
	Finally, $(Q_t\phi)(y):=q(t,y)\phi(y), (R_t\phi)(y):=r(t,y)\phi(y)$, and $(G\phi)(y):=g(y)\phi(y)$.  Theorems~\ref{Th4.4} and ~\ref{EnergyEquality} can be applied to solve the above quadratic optimal control  of the Anderson model.
\end{example} 

\begin{example} Consider the controlled Anderson system with the coefficient $b$ being time-invariant and $b^{-1}$ existing and being bounded. Then the system~\eqref{control system} is stablizable by the feedback control $u=-\lambda b^{-1}(y)X$ for sufficiently large  $\lambda $. To show this, we have for $\widetilde X_t:=e^{\lambda t} X_t$, 
	\begin{equation}
	{\widetilde	X}_t=e^{At}x+\int_0^t\sum_{j=1}^\infty e^{A(t-s)}  C_j{\widetilde X}_s\,d \beta_s^j.
	\end{equation}
	Therefore,  we have
	\begin{eqnarray}
	\mathbb{E}[|\widetilde{X}_t|^2]&\le&  C|x|^2+C\mathbb{E}\int_0^t\sum_{j=1}^\infty |e^{A(t-s)}  C_j{\widetilde X}_s|^2ds\nonumber\\
	&\le&  C|x|^2+C\int_0^t(t-s)^{-2\alpha}\mathbb{E}[|\widetilde{X}_s|^2]ds. 
	\end{eqnarray}
	From Gronwall's inequality, we have 
	$$
	\mathbb{E}[|\widetilde{X}_t|^2]\le  C |x|^2 e^{Ct}, \quad \mathbb{E}[|X_t|^2]\le  C |x|^2 e^{(C-\lambda)t}. 
	$$
	Therefore, Theorems~\ref{alg} and ~\ref{alg verification} can be applied to the Anderson model. 
	
\end{example}

\end{document}